\documentclass[a4paper,11pt]{article}

\usepackage{amsmath,latexsym,amssymb}

\newtheorem{thm}{Theorem}
\newtheorem{prop}[thm]{Proposition}
\newtheorem{lemma}{Lemma}
\newtheorem{cor}{Corollary}
\newtheorem{main}{Main Theorem}

\newcommand{\proof}[1][]{{\it Proof{#1}: }}
\newcommand{\qed}[1][1cm]{\hspace*{\fill} $\Box$ \vspace{#1}}

\renewcommand{\P}{{{\mathbf P}^1}}
\newcommand{\ZZ}{{\mathbf Z}}
\newcommand{\CC}{{\mathbf C}}

\newcommand{\LL}{{\mathbf L}}

\newcommand{\fk}{{\mathbf F}_k} % hirzebruch surface
\newcommand{\hirz}{{\mathbf F}}
\newcommand{\ph}{PH^0}
\newcommand{\isec}{C_{-k}} % section at infinity
\newcommand{\low}{^{\phantom{1}}}

\newcommand{\del}{\partial}
\newcommand{\ON}{\operatorname}

\renewcommand{\b}{\beta}

\newcommand{\e}{\varepsilon}

\renewcommand{\l}{\lambda}

\newcommand{\s}{\sigma}
\newcommand{\x}{\xi}

\newcommand{\tto}{\longrightarrow}
\newcommand{\inj}{\hookrightarrow}
\newcommand{\inv}{^{^{-1}}}
\newcommand{\surj}{\to\!\!\!\!\!\!\!\!\tto}

\newcommand{\bfami}{{\cal B}}

\newcommand{\dfami}{{\cal D}}

\newcommand{\ffami}{{\cal F}}

\newcommand{\ofami}{{\cal O}}

\newcommand{\wfami}{{\cal W}}
\newcommand{\xfami}{{\cal X}}

\newcommand{\labell}[1]{\label{#1}}

\newcommand{\br}{\ON{Br}}

\newcommand{\map}{M^0}
\newcommand{\psl}{\ON{PSL_2\ZZ}}
\newcommand{\slz}{\ON{SL_2\ZZ}}

\newcommand{\diff}{\ON{Diff}}

\newcommand{\pnull}{\l}
\newcommand{\qnull}{\xi}
\newcommand{\peins}{\l_1}
\newcommand{\qeins}{\xi_1}
\newcommand{\imag}{\sqrt{-1}}
\newcommand{\bif}{\ON{Bif}}
\newcommand{\bdiff}{\varphi} % diffeomorphism of the punctured base
\newcommand{\Phisl}{{\mit\Phi}}
\renewcommand{\Phisl}{{\tilde{\varphi}}}

\begin{document}

{
\Large\bf
\noindent
On bifurcation braid monodromy of elliptic fibrations\\[1.2cm]
}

{\bf
\noindent
Michael\ L\"onne\footnote{Institut f\"ur Mathematik, Universit\"at Hannover, Am
Welfengarten 1, 30167 Hannover, Germany (e-mail: loenne@math.uni-hannover.de)} 
}

\paragraph{Abstract.}
We define a monodromy homomorphism for irreducible families of regular
elliptic fibrations which takes values in the mapping class group of a punctured
sphere. In the computation we consider only elliptic fibrations which contain no
singular fibres of types other than $I_1$ and $I_0^*$.

We compare the maximal groups, which can be the monodromy groups of algebraic,
resp.\ differentiable families of elliptic surfaces, and give an algebraic
criterion for the equality of both groups which we can check to apply in case
the number of $I_1$ fibres is at most $6$.

\section*{introduction}

The monodromy problems we want to discuss fit quite nicely into the following
general scheme:
Given an algebraic object $X$ consider an algebraic family $g:\xfami\to T$ such
that a fibre $g\inv(t_0)$ is isomorphic to $X$ and such that the
restriction to
a connected subfamily $g|:\xfami'\to T'$ containing $X$ is a locally trivial
differentiable fibre bundle. If $G$ is the structure group of this bundle, the
geometric monodromy is the natural homomorphism $\rho:\pi_1(T',t_0)\to G$. A
monodromy map with values in a group $A$ is obtained by composition with some
representation $G\to A$.\\

In the standard setting $X$ is a complex manifold, e.g.\ a smooth complex
projective curve. In this case $\xfami$ is a flat family of compact curves
containing
$X$, the subfamily $\xfami'$ contains only the smooth curves and is a locally
trivial bundle of Riemann surfaces with structure group the
mapping class group $M\!\:\!ap(X)$. From the geometric monodromy one can
obtain the
algebraic monodromy by means of the natural representation
$M\!\:\!ap(X)\to Aut(H_1(X))$.\\

In the present paper we investigate the monodromy of regular elliptic
fibrations.
So $X$ is an elliptic surface with a map $f:X\to\P$ onto the projective line.
We consider families $g:\xfami\to T$ of elliptic surface containing $X$
with a map
$f_T:\xfami\to\P$ which extends $f$ and induces an elliptic fibration on each
surface. Subfamilies $\xfami'$ are to be chosen as differentiable fibre bundles
with
structure group $\diff_f(X)$, the group of isotopy classes of diffeomorphism
which commute with the fibration map up to a diffeomorphism of the base.

In the given setting there is a natural representation of $\diff_f(X)$ taking
values in the mapping class group of the base $\P$ punctured at the singular
values, \cite{birman},
\begin{eqnarray*}
\map_n & = & \left\langle \s_1,...,\s_{n-1}\,\left|\,
\begin{array}{l}
\s_i\s_{i+1}\s_i=\s_{i+1}\s_i\s_{i+1},\s_i\s_j=\s_j\s_i,|i-j|\geq2,\\[1.8mm]
\s_1\cdots\s_{n-2}\s_{n-1}^2\s_{n-2}\cdots\s_1=1=
(\s_1\cdots\s_{n-1})^n
\end{array}
\right.\right\rangle.
\end{eqnarray*}
Since there is a natural surjective homomorphism $\pi:\br_n\to\map_n$ from the
braid group on $n$ strands, we call the associated homomorphism the 
\emph{braid monodromy} of the family $\xfami'$. 

Essentially we consider generic fibrations with singular fibres of type $I_1$
only. To give some flavour of the general case we allow singular fibres also of
type $I_0^*$, cf.\ Kodaira's list \cite{ko}, \cite[p.150]{bpv}:

\setlength{\unitlength}{3.4mm}
\begin{picture}(40,10)
\linethickness{1pt}

%arcs
\bezier{200}(3,7)(10,0)(10,5)             
\bezier{200}(3,3)(10,10)(10,5)             

\bezier{200}(20,3)(22,5)(25,8)

\bezier{200}(21,2)(24,5)(26,7)

\bezier{200}(23,8)(27,5)(31,2)

\bezier{200}(27,3)(29,5)(32,8)

\bezier{200}(28,2)(31,5)(33,7)

\end{picture}
\vspace*{-9mm}
$$
I_1\hspace*{63mm}I_0^*\hspace*{23mm}
$$

We call a subgroup $E$ of a spherical mapping class group the {\sl braid
monodromy group}
of a fibration $X$, if $E$ is the smallest subgroup (w.r.t.\ inclusion)
such that
for all admissible $\xfami$ the image of the braid monodromy is a subgroup
of $E$ up to conjugation.

\begin{main}
\labell{main}
The braid monodromy group of a regular elliptic fibration $X$
with no singular fibres except $6k$ fibres of type $I_1$ and $l$ fibres of type
$I^*_0$ is a subgroup of $\map_{6k+l}$ representing the conjugacy class of
$$
\overline{E}_{6k,l}:=\left\langle \s_{ij}^{m_{ij}},\,i<j\,\left|\quad m_{ij}=\left\{
\begin{array}{ll}
1 & \text{ if }i,j\leq l \,\,\vee\,\, i\equiv j\,(2),i,j>l\\
2 & \text{ if } i\leq l<j\\
3 & \text{ if } i,j>l, i\not\equiv j \,(2)
\end{array}
\right.\right.\right\rangle.
$$
(Here $\s_{i,i+1}:=\s_i$, while for $j>i+1$ we define
$\s_{i,j}:=\s_{j-1}\cdots\s_{i+1}\s_i\s_{i+1}\inv\cdots\s_{j-1}\inv$.)
\end{main}

\newcommand{\ebar}{\overline{E}}
\newcommand{\kbar}{\overline{K}}

Of course it suffices to take an admissible family $\xfami$
which is topologically versal in the sense, that each admissible family is
topologically equivalent to a family induced from $\xfami$.
\\
\\
Our interest in this kind of results stems from three sources:
\\

Each family of fibrations we consider is a subfamily of a family of smooth
elliptic surfaces. Hence 
the natural surjection for the fundamental groups of the base spaces induces a
map from  $\pi_1^{\text{\sl surf}.}$ of the base of elliptic surfaces
to a quotient of $\ebar$:
$$
\begin{array}{ccccc}
K & \inj & \pi_1^{\text{\sl fib.}} & \surj & \pi_1^{\text{\sl surf}.}\\
\downarrow & & \downarrow & & \downarrow\\[1mm]
\kbar & \inj & \ebar & \surj & \ebar/ \kbar
\end{array}
$$
This diagram should be of some help for the understanding of the
fundamental group of some suitably defined moduli spaces of regular elliptic
surfaces.
\\

There is some analogy to the proposal of Donaldson \cite{d} to investigate
$\pi_1^{\text{\sl surf}.}$ by a suitable monodromy map to the
symplectic isotopy classes of symplectomorphisms.

In fact $\ebar/\kbar$ is a quotient of the fundamental group of a discriminant
complement for the semi-universal unfolding of a hypersurface singularity
$y^3+x^{6k}$, cf.\ \cite{habil}. This observation and the results
of Seidel, \cite{sei}, and Auroux, Munoz, Presas, \cite{amp} point at a
striking resemblance of the two homomorphism which deserves further
investigations.
\\

Finally we would like to know about the implications of our result in the
category of differentiable bundles of elliptic fibrations:
\\

Along with the proof of the theorem we will notice that each mapping class in 
the braid monodromy group $\overline{E}(X)$ is represented by a diffeomorphism
which can be lifted to a diffeomorphism of $X$ inducing the trivial mapping
class on some generic fibre. Hence we ask for the converse:
\begin{quote}
Does every diffeomorphism of $X$, isotopic to the identity mapping on some
generic fibre, induce a mapping class of the punctured base which is in the
monodromy group of $X$?
\end{quote}

A positive answer would yield a topological characterisation of the braid
monodromy group!\\

In fact we show that the group of mapping classes induced in the said way
is the image under $\pi$ of the stabiliser group $Stab_\psi$ 
of the algebraic monodromy $\psi$ with respect to an appropriate Hurwitz action.

Then we use \cite{hurwitz} to give an affirmative answer to the question above
in case the number of fibres of type $I_1$ does not exceed $6$.

\begin{thm}
\label{1}
Let $\overline{E}:=\overline{E}_{6,\ell}$ be the braid monodromy group of a
regular elliptic fibration $X$ with no singular fibres except $6$ fibres of
type $I_1$ and $\ell$ fibres of type $I_0^*$.\\
Then 
\begin{eqnarray*}
\overline{E} & = & \pi(\text{Stab}_\psi).
\end{eqnarray*}
\end{thm}

\section*{bifurcation braid monodromy}

With each locally trivial bundle one can associate the structure homomorphism
defined on the fundamental group of the base with respect to any base point. It
takes values in the mapping classes of the fibre over the base point.

Given a curve $C$ in the affine plane we can take a projection to the
affine line
which restricts to a finite covering $C\to\CC$. The complement of the curve
and its
vertical tangents is the total space of a punctured disc bundle over the
complement of the branch points in the affine line.

The structure homomorphism of this bundle is called the braid monodromy of the
plane curve with respect to the projection, and it can be naturally
regarded as a
homomorphism from the fundamental group of the branch point complement to the
braid group, since the latter is naturally isomorphic to the mapping class
group of the punctured disc.

This definition is readily generalized to the case of a divisor in the
Cartesian product of the affine or projective line with an irreducible base
$T$.
Then the structure homomorphism takes values in a braid group, resp.\ in a
mapping
class group of a punctured sphere $\map_n$.

In situations as we are interested in, such a divisor is defined as the
locus of
critical values of a family of algebraic functions of
constant bifurcation degree
with values in $\LL\cong\P$ or $\CC$.
Thus we give the relevant definitions:

\begin{description}

\item[Definition:]
A flat family $\xfami\to T$ with an algebraic morphism $f:\xfami\to \LL$ is
called a
{\it framed family of functions} $(\xfami,T,f)$.

\item[Definition:]
The {\it bifurcation set} of a framed family of functions
over $T$ is the smallest Zariski closed subset $\bfami$ in $T\times \LL$ such
that the diagonal map $\xfami\to T\times \LL$ is smooth over the complement of
$\bfami$.

\item[Definition:]
The {\it discriminant set} of a framed family of functions over $T$ is the
divisor in $T$ such that the bifurcation set $\bfami$ is an unbranched cover
over its complement by the restriction of the natural projection $T\times \LL\to
T$.

\item[Definition:]
A framed family of functions is called {\it of constant
bifurcation degree} if the bifurcation set is a finite cover of $T$.

\item[Definition:]
The {\it bifurcation braid monodromy} of a framed family of functions with
constant bifurcation degree over an irreducible base $T$ is defined to be the
braid monodromy of $\bfami$ in $T\times\LL$ over $T$.

\end{description}

Note that this definition of braid monodromy differs slightly from the
definition
given in the introduction but that the resulting objects are the same.

\section*{families of divisors in Hirzebruch surfaces}

Given a Hirzebruch surface $\fk$ with a unique section $\isec$ of
selfintersection $-k$, we consider families of divisors on which the ruling of
the Hirzebruch surface defines families of functions with constant bifurcation
type.\\
We can pull back divisors form the base along the ruling to get divisors on
$\fk$ which we call vertical, among others the fibre divisor $L$.\\
Consider now the family of divisors on $\fk$ which consist of a vertical part
of degree $l$ and a divisor in the complete linear system of
$\ofami_{\fk}(4\isec+3kL)$ called the horizontal part. It is a family
parameterized by $T=\ph(\ofami_{\P}(l))\times\ph(\ofami_{\fk}(4\isec+3kL))$
with
total space
$$
\dfami_{k,l}=\left\{(x,t)\in\fk\times T\,|\,x\in D_t\subset\fk\right\}.
$$
Let $T'$ be the Zariski open subset of $T$ which is the base of the family
$\dfami'_{k,l}$ of divisors in $\dfami_{k,l}$ with reduced horizontal part.

\begin{lemma}
\labell{dfami}
The ruling on $\fk$ defines a morphisms $\dfami_{k,l}'\to\P$ by which it
becomes a framed family of functions of constant bifurcation degree.
\end{lemma}

\proof
The critical value set of the vertical part of a divisor is the divisor of
which it is the pull back, thus it is constant of degree $l$.\\
The assumption on reducedness forces the horizontal part to be without fibre
components. We may even conclude that a reduced horizontal part consists of
$\isec$ and a disjoint divisor which is a branched cover of the base of degree
$3$. The critical values set is therefore the branch set which is of constant
degree $6k$, and we are done.
\qed

The abstract braid group given by the presentation
\begin{eqnarray*}
\br_n & = & \langle \s_1,...,\s_{n-1}\,|\, \s_i\s_{i+1}\s_i=\s_{i+1}\s_i\s_{i+1},
\s_i\s_j=\s_j\s_i \text{ if } |i-j|\geq2\rangle.
\end{eqnarray*}
and the abstract mapping class group $\map_n$ given by the presentation
in the introduction are naturally identified with the mapping
class group of the punctured disc, resp.\ sphere. Such an identification is
given if each
$\s_i$ is realised by the half-twist on an embedded arc $a_i$ connecting two
punctures provided that
$a_i\cap a_{i+1}$ is a single puncture and $a_i\cap a_j$ is empty if
$|i-j|\geq2$.

\begin{prop}
\labell{dmono}
The image of the bifurcation braid monodromy homomorphism of the family
$\dfami'_{k,l}$ is conjugated to the subgroup of $\map_{6k+l}$:
$$
\overline{E}_{6k,l}:=\left\langle \s_{ij}^{m_{ij}},\,i<j\,\left|\quad m_{ij}=\left\{
\begin{array}{ll}
1 & \text{ if }i,j\leq l \,\,\vee\,\, i\equiv j\,(2),i,j>l\\
2 & \text{ if } i\leq l<j\\
3 & \text{ if } i,j>l, i\not\equiv j \,(2)
\end{array}
\right.\right.\right\rangle.
$$
\end{prop}

The proof of this proposition and a couple of preparatory results will take the
rest of the section.\\

First note that our whole concern lies in the understanding of the bifurcation
set $\bfami$ in $T'\times\P$ with its projection to $T'$. As an approximation
we will consider families of affine plane curves given by families of
polynomials in affine coordinates $x,y$ with the regular map induced by the
affine projection $(x,y)\mapsto x$.

Their bifurcation sets are contained in the Cartesian product of the family
bases with the affine line $\CC$, and it will soon be shown that this pair can
be induced from $(T'\times\P,\bfami)$.
Eventually we can extract all necessary information from such families to
prove our claim.

\begin{lemma}
\labell{div/discr}
Consider $y^3-3\,p(x)\,y\,+\,2\,q(x)$ as a family of polynomial functions
\mbox{$\CC^2\times T\to\CC$} parametrised by a base $T$ of pairs $p,q$ of
univariate polynomials. Then the bifurcation set is the zero set of
$g(x):=p^3(x)-q^2(x)$, the discriminant set is the zero set of the discriminant
of $g$ with respect to
$x$.
\end{lemma}

\proof
The bifurcation divisor is cut out by the discriminant polynomial of
$y^3-3p(x)y+2q(x)$ with respect to $y$. The first claim is then immediate since
$g$ is proportional to the corresponding Sylvester determinant:
$$
\begin{vmatrix}
1 & 0 & -3p & 2q \\
 & 1 & 0 & -3p & 2q \\
3 & 0 & -3p\\
 & 3 & 0 & -3p\\
 & & 3 & 0 & -3p
\end{vmatrix}
$$
For the second claim we note that a pair $p,q$ belongs to the 
discriminant set if and only if $p^3-q^2$ has a
multiple root hence this locus is cut out by the discriminant of
$g$ with respect to $x$.
\qed

\begin{lemma}
\labell{dis-comp}
The discriminant locus of a family $y^3+3r(x)y^2-3p(x)y+2q(x)$ is the union of
the degeneration component of triples $p,q,r$ defining singular curves and the
cuspidal component of triples defining polynomial maps with a degenerate critical
point.
\end{lemma}

\proof
In general a branched cover of $\CC$ has not the maximal number of branch points
only if the cover is singular, or the number of preimages of a branch point
differs by more than one from the degree of the branching. The second
alternative
occurs only if there is a degenerate critical point in the preimage or if
there are
two critical points. Since the last case can not occur in a cover of degree
only
three we are done.
\qed

\begin{lemma}
\labell{cusp-mult}
Given the family $y^3+3r(x)y^2-3p(x)y+2q(x)$ the cuspidal component of the
discriminant is the zero set of the resultant of $p(x)+r^2(x)$ and
$2q(x)-r^3(x)$ with respect to $x$.\\ 
Its equation - considered a polynomial in the variable $\l_0$ - is of degree $n$
with coprime coefficients if\\[-2mm]
$$
p(x)=\sum_{i=0}^d\l_ix^i,\quad q(x)=x^n+\sum_{i=0}^{n-1}\x_ix^i,\quad
r(x)=\sum_{i=0}^{\lfloor n/3\rfloor}\zeta_ix^i.
$$
\end{lemma}

\proof
The cuspidal discriminant is the locus of all parameters for which there is a
common zero of $f,\del_yf,\del^2_yf$. Since $\del^2_yf=0$ is linear in $y$,
we can eliminate $y$ and get the resultant of $p(x)+r^2(x)$ and $2q(x)-r^3(x)$
with respect to $x$.\\
By the degree bound on $q$ and $r$ the discriminant equation is the resultant
of a matrix in which the variable $\l_0$ occurs exactly $n$ times. Moreover
the diagonal determines the coefficient of $\l_0^n$ to be
a power of 
$(1-\zeta_{n/3}^3)$ resp.\
$1$ depending on whether $n/3\in\ZZ$ or not. Even in the first case the
coefficients are coprime since the resultant is not divisible by
$(1-\zeta_{n/3}^3)$. 
\qed

\begin{lemma}
\labell{dgn-mult}
For the family
$y^3+3r(x)y^2-3p(x)y+2q(x)$ 
the degeneration component of the discriminant is
the locus of triples for which there is a common zero in
$x,y$ of the polynomial and its two partial derivatives.
Its equation - considered a polynomial in the variable $\xi_0$ - is monic of
degree $2n-2$ if\\[-4mm]
$$
p(x)=\sum_{i=0}^d\l_ix^i,\quad q(x)=x^n+\sum_{i=0}^{n-1}\x_ix^i,\quad
r(x)=\sum_{i=0}^{\lfloor n/3\rfloor}\zeta_ix^i.
$$
\end{lemma}
\proof
The degeneration locus is given by the Jacobian criterion as claimed. The
equation of the discriminant of the subdiagonal unfolding of the
quasi-homogeneous singularity $y^3-x^n$ is known to the quasi-homogeneous and of
degree $2n-2$ in $\xi_0$. Since the unfolding over the $\xi_0$-parameter is a
Morsification the coefficient of $\xi^{2n-2}_0$ must be constant.
\qed

\begin{lemma}
\labell{special}
The bifurcation braid monodromy of the family $y^3-3\pnull y+2(x^n+\qnull)$ maps
onto a subgroup of $\br_{2n}$ which is conjugated to the subgroup generated by
$$
(\s_1\cdots \s_{2n-1})^3,(\s_{2n-2,2n}\cdots \s_{2,4})^{n+1},
(\s_{2n-3,2n-1}\cdots \s_{1,3})^{n+1}.
$$
\end{lemma}

\proof
As one can show with the help of the preceding lemmas, the discriminant
locus is the union of the degeneration locus and the cuspidal component
which are
cut out respectively by the polynomials $\pnull^3-\qnull^2$ and $\pnull$.

By Zariski/van Kampen the fundamental group of the complement with base point
$(\pnull,\qnull)=(1,0)$ is generated by the fundamental group of the complement
restricted to the line $\pnull=1$ and the homotopy class of a loop which links
the line $\pnull=0$ once.

For the pair $(1,0)$ the set of regular values of the polynomial consists
of the
affine line punctured at the $(2n)$-th roots of unity, which we number
counterclockwise, $1$ the first puncture.
To express the bifurcation braid group in terms of abstract generators, we
identify the elements $\s_{i}$ with the half twist on the circle segment
between the $i$-th and $i+1$-st puncture.

For the line $\pnull=1$ the bifurcation locus is given by
$(x^n+\qnull-1)(x^n+\qnull+1)$. This locus is smooth but branches of degree $n$
over the base at $\qnull=\pm1$. The corresponding monodromy transformations are
the second and third transformation given in the claim.

Associated to the degeneration path $(\pnull,\qnull)=(1-t,\imag t)$, $t\in[0,1]$
there is a loop in the complex line $\pnull=1+\imag\qnull$ which links the line
$\pnull=0$. For this degeneration the
bifurcation divisor is regular and contains points of common absolute value
determined by $t$ only, except for $t=1$ where it has $n$ ordinary cusps with
horizontal tangent cone. Since a cusp corresponds to a triple half twist
and the first and second puncture merge in the degeneration, the monodromy
transformation for our loop is the first braid of the claim.
\qed

\begin{lemma}
\labell{semi-gen}
The bifurcation braid monodromy of the family $y^3-3\pnull y+2(x^n+\qnull+\e x)$,
$\e$~small and
fix, is in the conjugation class of the subgroup of the braid group
generated by
$$
(\s_{1}\s_3\cdots \s_{2n-1})^3,\s_{i,i+2},i=1,...,2n-2.
$$
\end{lemma}

\proof
The discriminant locus in the $\pnull,\qnull$ parameter plane consists again of the
cuspidal component $\pnull=0$ and the degeneration component.

Since the perturbation $\e$ is arbitrarily small, some features of the
family of
lemma~\ref{special} are preserved. The conclusion of the Zariski/van Kampen
argument still holds,
each braid group generators $\s_{i}$ is now realized as half twist on
segments of a slightly distorted circle, and the loop linking $\pnull=0$ is only
slightly perturbed. So the monodromy transformation associated to this loop is
formally the same as before, the first braid in the claim.

The dramatic change occurs in the bifurcation curve over the line $\pnull=1$.
Now the
bifurcation locus is the union of two disjoint smooth components each of which
branches simply of degree $n$ with all branch points near $\qnull=1$, resp.\
$\qnull=-1$.
Since the local model $x^n+\e x$ has the full braid group as its monodromy
group,
the monodromy along $\pnull=1$ is generated by the elements $\s_{i,i+2}$ as claimed.
\qed

\begin{lemma}
\labell{generic}
The bifurcation braid monodromy of the
family
$$
y^3-3(\pnull+\peins x)y+2(x^n+\qnull+\qeins x)
$$
is in the conjugation class of the subgroup generated by
$$
\s_{i}^{3},i\equiv1(2),\s_{i,i+2},0<i<2n-1.
$$
\end{lemma}

\proof
The components of the discriminant are the degeneration component and
the cuspidal component.

The line $\pnull=1,\peins=0,\qeins=\e$ small and fix, is generic for the
degeneration component
and we may conclude from lemma \ref{semi-gen} that there are elements in the
fundamental group of the discriminant complement with respect to
$(\pnull,\peins,\qnull,\qeins)=(1,0,0,0)$ which map to $\s_{i,i+2}$ as in lemma
\ref{semi-gen}.

Since the line $\qnull=i,\qeins=\peins=0$ is transversal for the cuspidal component
so are
parallel lines with $\peins=\e'$ small and fix. The bifurcation curve is then
given by
$(\pnull+\e' x)^3=(x^n+\imag)^2$. For $\pnull=0$ the critical values are
distributed in pairs along a circle in the affine line which merge pairwise for
$\peins\to0$.

But for $\e'$ sufficiently small and by varying $\pnull$ to an appropriate
small extend such
that the degeneration component is not met, there are $n$ obvious degenerations
when $\pnull$ is of the same modulus as $\e'$ and $\pnull+\e' x$ is a
factor of
$x^n+\imag$. By the local nature of the degeneration a degree argument shows that
these are all possibilities. Moreover one can easily see that the corresponding
monodromy transformations are triple half twist for each one of the pairs.

Moreover these twists are transformed to the transformations $\s_{i,i+1}^3$, $i$
odd, when transported along $(1-t,\imag t,0,\e't)$, $t\in[0,1]$ to the chosen
reference point.

The monodromy group is thus completely determined since the fundamental
group is
generated by elements which map to the given group under the monodromy
homomorphism.
\qed

\begin{lemma}
\labell{big+gen}
Let a family of plane polynomials be given which is of the form
\begin{equation*}
\label{pol}
y^3+3\Big(\sum_{i=0}^{d_r}\zeta_ix^i\,\Big)y^2
-3\Big(\sum_{i=0}^{d_p}\l_ix^i\,\Big)y+
2\big(x^n+\sum_{i=0}^{n-1}\xi_ix^i\,\big),
\end{equation*}
\hfill$2\leq n,0<3d_p\leq 2n,3d_r\leq n$.\\
Then the bifurcation braid monodromy group
is in the conjugation class of the subgroup of the
braid group generated by
$$
\s_{i}^{3},i\equiv1(2),\s_{i,i+2},0<i<2n-1.
$$
\end{lemma}

\proof
Since the family considered in the previous lemma is a subfamily now and
has the
claimed monodromy, we have to show that the new family has no additional
monodromy
transformations.

In the proof above we have seen that the cuspidal component is cut in $n$
points
by a line in $\l_0$ direction. The component is reduced since its
multiplicity at the origin is $n$, too, by lemma \ref{cusp-mult}. The
degeneration component is reduced by the analogous argument relying on lemma
\ref{dgn-mult} and transversally cut in $2n-2$ points by a line in $\xi_0$
direction.
Hence by Zariski/van Kampen arguments as proved in \cite{bessis} the fundamental
group of the discriminant complement of the subfamily surjects onto the
fundamental group of the family considered now.
\qed[4mm]

\begin{lemma}
\labell{red-gen}
Define a subgroup of the braid group $\br_{2n+l}$:
$$
E_{2n,l}:=\left\langle \s_{i,j}^{m_{ij}},\,i<j\,\left|\quad m_{ij}=\left\{
\begin{array}{ll}
1 & \text{ if }i,j\leq l \vee i\equiv j (2),i>l\\
2 & \text{ if } i\leq l<j\\
3 & \text{ if } i,j>l \wedge i\not\equiv j (2)
\end{array}
\right.\right.\right\rangle.
$$
Then the same subgroup is generated also by
$$
\s_{i},i<l,\s_{i,i+2},l<i,\s^2_{i,j},i\leq l<j,\s^3_{i},i>l,i\not\equiv
l (2).
$$
\end{lemma}

\proof
We have to show that the redundant elements can be expressed in the
elements of the
bottom line. This is immediate from the following relations $(i<j)$:
$$
\begin{array}{lll}
\s_{i,j} & = & \s\inv_{j-1}\cdots \s\inv_{i+1} \s_{i}\low
\s_{i+1}\low\cdots \s_{j-1}\low
,\quad j\leq l,\\
\s_{i,j} & = & \s\inv_{j-2,j}\cdots \s\inv_{i+2,i+4} \s_{i,i+2}\low
\s_{i+2,i+4}\low\cdots \s_{j-2,j}\low
,\quad l< i,i\equiv j(2),\\
\s^3_{i,j} & = & \s\inv_{j-2,2}\cdots \s\inv_{i+1,i+3} \s^3_{i}
\s_{i+1,i+3}\low\cdots
\s_{j-2,j}\low ,\quad l<i,i\not\equiv l,j(2),\\
\s^3_{i,j} & = & \s\inv_{j-2,2}\cdots \s\inv_{i+1,i+3} \s_{i-1,i+1}\low
\s^3_{i-1,i}
\s\inv_{i-1,i+1} \s_{i+1,i+3}\low\cdots \s_{j-2,j}\low
 ,\quad l<i,j\not\equiv l,i(2).
\end{array}
$$
\qed

\begin{lemma}
\labell{mono-full}
Consider the family $\left(y^3-3p(x)y+2q(x)\right)a(x)$
parametrised by triples $p,q,a$, with $p$ from the vector space of univariate
polynomials of degree at most $2n/3$, $q,a$ from the affine space of monic
polynomials of degree $n$ and $l$ respectively.
Then the subgroup $E_{2n,l}$ of $\br_{2n+l}$ is conjugate to a subgroup of the
image of the bifurcation braid monodromy.
\end{lemma}

\proof
We choose our reference divisor to be $(y^3-3y+2x^n)\prod_i^l(x-l-2+i)$ with
corresponding bifurcation set $x_i=l+2-i,i\leq l$ on the real axis and
$x_{l+1}=1$ and
the $x_i,i>l+1$ equal to the $2k^{\text{th}}$-roots of unity in
counterclockwise
numbering.\\
We identify the elements $\s_{i,j}$ of the braid group with the half twist
on arcs
between $x_i,x_j$, which are chosen to be
\begin{enumerate}
\item
a circle segment through the lower half plane, if $i,j\leq l$,
\item
a circle secant in the unit disc, if $i,j>l$,
\item
the union of a secant in the unit disc to a point on its boundary between
$x_{2n+l}$
and $1$ with an arc through the lower half plane, if $i\leq l<j$.
\end{enumerate}
(Of each kind we have depicted one in the following figure.)

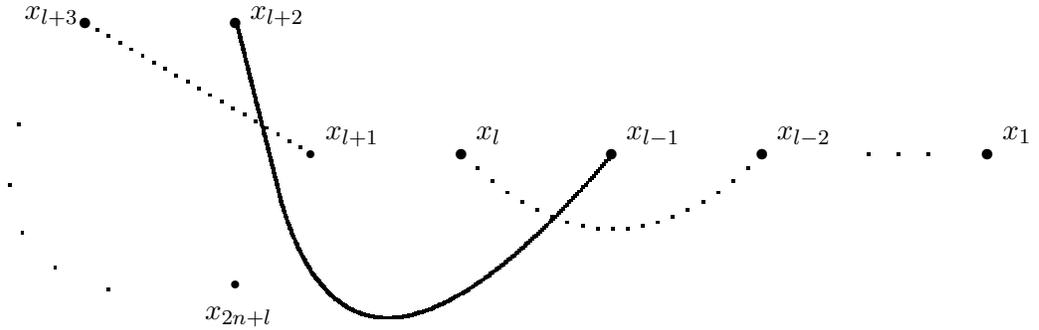
\setlength{\unitlength}{2mm}
\begin{picture}(68,22)
\linethickness{1pt}
%points
\put(15,1.3){\circle*{.6}}
%\put(10,10){\circle*{.7}}
\put(20,10){\circle*{.5}}
\put(30,10){\circle*{.7}}
\put(40,10){\circle*{.7}}
\put(50,10){\circle*{.7}}
\put(65,10){\circle*{.7}}
\put(5,18.7){\circle*{.7}}
\put(15,18.7){\circle*{.7}}

%dots
\bezier{4}(.5,12)(-2,3)(6.5,1)
\bezier{2}(57,10)(59,10)(61,10)

%arcs
\bezier{20}(30,10)(40,0)(50,10)              % type 1
\bezier{20}(20,10)(12.5,14.35)(5,18.7)       % type 2

\bezier{120}(15,18.7)(16.5,12.7)(18,6.7)     % type 3
\bezier{240}(18,6.7)(23,-10)(40,10)

%labels
%\put(10,10){0}
\put(51,11){$x_{l-2}$}
\put(66,11){$x_1$}
\put(16,19){$x_{l+2}$}
\put(1,19){$x_{l+3}$}
\put(21,11){$x_{l+1}$}
\put(31,11){$x_l$}
\put(41,11){$x_{l-1}$}
\put(13,-1){$x_{2n+l}$}

\end{picture}
\vspace*{6mm}\\

Since keeping the horizontal part $y^3-3y+2x^n$ fix, the
bifurcation divisor of the vertical is equivalent to that of the universal
unfolding of the function $x^l$ we have the elements $\s_{i},i<l$ in the
braid monodromy.
These elements are obtained for example in families
$$
a(x)=\left((x-l+i-3/2)^2+\l\right)\prod_{j\neq i,i+1}^l(x-l-2+j).
$$

The second set of elements, $\s_{i,j}^2,i\leq l<j$ is obtained by families
of the
kind
$$
(y^3-3y+2x^n)(x-l-2+i-\l)\prod_{j\neq i}^l(x-l-2+j)
$$
since the zero $l+2-i+\l$ may trace any given path in the range of the
projection, in particular that around an arc on which the full twist
$\s_{i,j}^2$
is performed.\\
Finally varying the horizontal part as in lemma \ref{generic}
while keeping the
$a(x)$ factor fix proves that the braid group elements $\s_{i,i+2},l<i$ and
$\s_{i},l<i,i\not\equiv l(2)$ are in the image of the monodromy. So we may
conclude that this image contains $E_{2n,l}$ up to conjugacy.
\qed

\proof[ of prop.\ \ref{dmono}]
Denote by $S$ the Zariski open subset of $T'$ which parameterizes divisors
of the
family $\dfami_{k,l}'$ which have no singular value at a point $\infty\in\P$.
The corresponding family in $\fk\times S$ may then be restricted to a
family $\ffami_{k,l}$ in $\CC\times\CC\times S$, where $\fk$ is
trivialized as $\CC\times\CC$ in the complement of the negative section
$\isec$ and the fibre over $\infty$.
By construction $\ffami_{k,l}$ has constant bifurcation degree.\\
Consider now the family of polynomials
$$
\left( y^3+3r(x)y^2-3p(x)y+2q(x)\right)a(x),
$$
where $r,p,q,a$ are taken from the family of all quadruples of polynomials in
one variable subject to the conditions that
\begin{enumerate}
\item
$r,p$ are of respective degrees $k$ and $2k$,
\item
$q$ is monic of degree $3k$, $a$ is monic of degree $l$
\item
the discriminant of $y^3+3r(x)y^2-3p(x)y+2q(x)$ is not identically zero.
\end{enumerate}
This family can be naturally identified with $\ffami_{k,l}$.
By lemma \ref{mono-full}, up to conjugacy, $E_{6k,l}$ is contained in the
monodromy image $\rho(\pi_1(S\setminus Discr(\ffami_{k,l})))$.

For the converse we note that the bifurcation set
of the family decomposes into the bifurcation sets $\bif_h$ of the horizontal
part $y^3+3r(x)y^2-3p(x)y+2q(x)$ and $\bif_v$ of the vertical part $a(x)$. Hence
the monodromy is contained in the subgroup $\br_{(6k,l)}$ of braids which do not
permute points belonging to different components. $\br_{(6k,l)}$ has natural
maps to $\br_{6k}$ and $\br_l$ which commute with the braid monodromies of both
bifurcation set considered on their own.

The discriminant decomposes into the discriminants of $\bif_h$,
$\bif_v$ and the divisor of parameters for which $\bif_h\cap\bif_v$ not empty.
They give rise in turn to braids which can be considered as elements in
$$
\br_{6k},\,\br_l\text{ resp.\ }\br_{(6k,l)}^{0,0}:=\{\b\in\br_{(6k,l)}|\,
\b\text{ trivial in }\br_{6k}\times\br_l\}.
$$
Now with lemma \ref{red-gen} we can identify $E_{6k.l}$ as
the subgroup of $\br_{6k+l}$ generated by
$E_{6k}\subset\br_{6k}$, $\br_l$ and $\br_{(6k,l)}^{0,0}$
which are generated in turn by the elements 
$$
\{\s_{i,i+2},\s^3_{i},l<i\},\{\s_{i},i<l\},\{\s^2_{i,j},i\leq l<j\}\text{
resp.}
$$
And by lemma \ref{big+gen} the image can not contain more elements.\\

Since the bifurcation diagram of $\ffami_{k,l}$ embeds in the
bifurcation diagram of $\dfami_{k,l}'$ with complement of codimension one, there
is a commutative diagram
$$
\begin{array}{ccc}
\pi_1(S\setminus Discr(\ffami_{k,l})) & \tto\!\!\!\!\!\!\!\!\tto
& \pi_1(T'\setminus Discr(\dfami_{k,l}'))\\
\downarrow & & \\[-4.7mm]
\downarrow & & \Big|\\[-.1mm]
E_{6k,l} & & \big|\\[-.5mm]
\rceil & & \Big|\\[-3mm]
\hspace*{.2mm}\downarrow & & \hspace*{.1mm}\big\downarrow\\
\br_{6k+l} & \tto\!\!\!\!\!\!\!\!\tto & \map_{6k+l}
\end{array}
$$
from which we read off our claim.
\qed

\begin{cor}
For any element $\b$ in the braid monodromy group of $\dfami_{k,l}'$ there is a
diffeomorphism of the base $\P$ which fixes a neighbourhood of
$\infty\in\P$ and
which represents the mapping class $\b$.
\end{cor}

\proof
The element $\b$ is image of an element $\b'$ in the braid monodromy of the
bifurcation diagram of $\ffami_{k,l}$. The bifurcation set does not meet the
boundary so integration along a suitable vector field yields a realisation of
$\b'$ as a diffeomorphism acting trivially on a neighbourhood of the boundary.
Its trivial extension to the point $\infty$ is the diffeomorphism sought for.
\qed

\section*{families of elliptic surfaces}

In this section we start investigating families of regular elliptic surfaces
for which the type of singular fibres is restricted to $I_1$ and $I_0^*$.
We will go back and forth between a family of elliptic fibrations, its
associated
family of fibrations with a section and a corresponding Weierstrass model
of the
latter, so we note some of their properties:

\begin{prop}
\labell{jac}
Given a family of elliptic fibrations with constant bifurcation type over an
irreducible base $T$, there is a family of elliptic fibrations with a section,
such that the bifurcation sets of both families coincide.
\end{prop}

\proof
Given a family as claimed there is the associated family of Jacobian
fibrations,
cf.\ \cite[I.5.30]{fm}. The bifurcation sets of both families coincide.
\qed

In turn, for each family of elliptic fibrations with a section there is a
corresponding family of Weierstrass fibrations, cf.\ Miranda \cite{mir}.

A regular Weierstrass fibration $W$ is defined by an equation
$$
wz^2=4y^3-3Pw^2y+2Qw^3
$$
in the projectivisation of the vector bundle
$\ofami\oplus\ofami(2\chi)\oplus\ofami(3\chi)$  over the projective line $\P$
where $\chi$ is the holomorphic Euler number of the fibration,
$w,y,z$ are 'homogeneous coordinates' of the bundle, and $P,Q$ are sections of
$\ofami(4\chi),\ofami(6\chi)$ respectively.\\
So $W$ is a double cover of the Hirzebruch surface $\hirz_{2\chi}={\mathbf
P}(\ofami\oplus\ofami(2\chi))$ branched along the section $\s_{2\chi}$ and the
divisor in its complement $\ofami(2\chi)$ defined by the equation
$y^3-3Py+2Q=0$.\\

A framed family of Weierstrass fibrations over a parameter space $T$ is a
given by
data as before where now $P,Q$ are sections of the pull backs to $T\times\P$
of $\ofami(4\chi),\ofami(6\chi)$ such that for each parameter $\l\in T$ they
define a Weierstrass fibration. In the sequel $P,Q$ are referred to as the
coefficient data of the Weierstrass family.

\begin{lemma}
\labell{w-fact}
Let $\wfami$ be the Weierstrass family associated to a framed family over
$T$ of
regular elliptic fibrations in which all surfaces have no singular fibres
except
for $l$ of type $I_0^*$ and $6k$ of $I_1$ with coefficient data $P,Q$, then
there
are three families of sections $a,p,q$ of $\ofami(l),\ofami(2k),\ofami(3k)$
respectively, such that $p,q$ have no common zero,
\begin{eqnarray*}
& p\cdot a^2=P, & q\cdot a^3=Q,
\end{eqnarray*}
and the bifurcation set is given by
\begin{eqnarray*}
a\left(p^3-q^2\right)=0 & \subset & T\times\P.
\end{eqnarray*}
\end{lemma}

\proof
By the classification of Kas \cite{kas} at base points of regular fibres the
discriminant $P^3-Q^2$ does not vanish, at base points of fibres of type $I_1$
the discriminant vanishes but neither $P$ nor $Q$ and at base points of fibres
of type $I_0^*$  the vanishing order of $P$ is two, the vanishing order of
$Q$ is
three.\\
Since by hypothesis the locus of base points of singular fibres of type $I_0^*$
form a family of point divisors of degree $l$ there is a section $a$ of
$\ofami(l)$ such that $P$ has a factor $a^2$ and $Q$ a factor $a^3$.\\
With $deg P=2(l+k),
deg Q=3(l+k)$ we get the other degree claims.\\
Finally the discriminant of the Weierstrass fibration is given by $P^3-Q^2$
which
has -- by the above -- the same zero set as $a\left(p^3-q^2\right)$.
\qed

\begin{description}
\item[Remark:] In the situation of the lemma, a family of divisors is given for
$\hirz_{k}$ by the equation $a(y^3-3pw^2y+2qw^3)=0$, $a$ cutting out the
vertical part. The double cover along this divisor is a family of fibrations
obtained from the original family by contracting all smooth rational curves of
selfintersection $-2$, of which there are four for each fibre of type $I_0^*$.
\end{description}

We are now prepared to come back to the main theorem:\\

\proof[ of the main theorem]
Given any framed family of regular elliptic fibrations containing $X$ we
consider
a Weierstrass model $\wfami$ of the associated Jacobian family. Since $\wfami$
is again framed there is an induced family of divisors on a Hirzebruch
surface obtained as before.\\
This family of divisors is a pull back from the space $\dfami_{k,l}$ so the
monodromy is a subgroup of the bifurcation monodromy of Hirzebruch divisors.

On the other hand for the family of triples of polynomials
$p(x),q(x),a(x)$ with $p$ of degree at most $2k$ and $q,a$ monic of degree $3k$
respectively $l$, we can form the family given by
$$
z^2=y^3-3p(x)a^2(x)+2q(x)a^3(x),
$$
which is Weierstrass in the complement of parameters where
$a(x)\left(p^3(x)-q^2(x)\right)$ has a multiple root or vanishes identically.
At least after suitable base change, cf.\
\cite[p.\ 163]{fm}, this Weierstrass family has a simultaneous resolution
yielding
a family $\xfami_{k,l}$ of elliptic surfaces with a section.\\
The Jacobian of $X$ is contained in $\xfami_{k,l}$, since its Weierstrass data
consist of sections $P,Q$ which are factorisable as $a^2p,a^3q$ according to
lemma \ref{w-fact} and after the choice of a suitable $\infty$ this data
can be identified with polynomials in this family.\\
The fibration $X$ is deformation equivalent to its Jacobian with
constant local analytic type, cf.\ \cite[thm.\ I.5.13]{fm} and hence of
constant
fibre type.
The monodromy group therefore contains the bifurcation monodromy group of
divisors on Hirzebruch surfaces $\dfami_{k,l}$ and so the two groups even
coincide.
\qed

Regarding elements in the braid monodromy as mapping classes again they can be
shown to be induced by diffeomorphism of the elliptic fibration, but more
is true
in fact:

\begin{prop}
\labell{lift}
For each braid $\b$ in the framed braid monodromy group there is a
diffeomorphism
of the elliptic fibration which preserves the fibration, induces $\b$ on
the base
and the trivial mapping class on some fibre.
\end{prop}

\proof
As we have seen in the corollary to prop.\ \ref{dmono} we can find a
representative $\bdiff$ for the braid $\b$ by careful integration of a suitable
vector field such that $\bdiff$ is the identity next to a point $\infty$.\\
In \cite[II.1.2]{fm} there is a proof for families of nodal elliptic fibrations
and sufficient hints for more general families of constant singular fibre
types,
that a horizontal vector field on the total family can be found which fails
to be
a lift only in arbitrarily small neighbourhoods of singular points on singular
fibres. Integration of such a vector field yields a diffeomorphism $\Phisl$ which
is a lift of $\bdiff$.\\
We have seen that the monodromy generators arising from the horizontal part can
be realized over a suitable polydisc parameter space, cf.\ lemma \ref{generic}.
Since the vertical part as in lemma \ref{mono-full} does not have any effect on
the fibre $F_\infty$ over $\infty$ we can conclude that this fibration
family is
the trivial family next to $F_\infty$. So we apply the argument above to get a
lift $\Phisl$ which induces the trivial mapping class on $F_\infty$.
\qed[2mm]

\section*{Hurwitz stabilizer groups}

In this section we determine the stabilizers of the action of the braid group
$\br_n$ on homomorphisms defined on the free group $F_n$ generated by elements
$t_1,...,t_n$.
The action is given by precomposition with the Hurwitz automorphism of $F_n$
associated to a braid in $\br_n$:
$$
\br_n\to\ON{Aut} F_n:\quad
\s_{i,i+1}\mapsto\left(t_j\mapsto
\begin{cases}
t_j & j\neq i,i+1\\
t_it_jt_i\inv & j=i\\
t_i & j=i+1
\end{cases}
\right).
$$
We start with a result from \cite{hurwitz}:

\begin{prop}
\labell{l}
Let $F_n:=\langle t_i,1\leq i\leq n\,|\quad\rangle$ be the free group on $n$
generators, define a homomorphism $\phi_n:F_n\to\br_3=\langle
a,b\,|\,aba=bab\rangle$ by
$$
\phi_n(t_i)=\left\{
\begin{array}{ll}
a & i\text{ odd}\\
b & i\text{ even}
\end{array}\right.
$$
and let $\br_n$ act on homomorphisms $F_n\to\br_3$ by Hurwitz
automorphisms of $F_n$. Then the stabilizer group $Stab_{\phi_n}$
contains the braid subgroup
$$
E_n=\langle \s_{i,j}^{m_{ij}}\:|\: m_{ij}=1,3
\text{ if $j\equiv i$, resp.\ $i\not\equiv j \mod 2$}\rangle
$$
with $E_n=Stab_{\phi_n}$, if $n\leq6$.
\end{prop}

Note that the action in \cite{hurwitz} was defined on tuples
$\left(\phi_n(t_1),...,\phi_n(t_n)\right)$ but that it is obviously
equivalent to
the action considered here.\\
This result can now be applied to find stabilizers of similar homomorphisms:

\begin{prop}
\labell{slcor}
Let $F_n:=\langle t_i,1\leq i\leq n\,|\quad\rangle$ be the free group on $n$
generators, define a homomorphism $\psi_n:F_n\to\slz$ by
$$
\psi_n(t_i)=\left\{
\begin{array}{ll}
\left(\begin{smallmatrix}1&1\\0&1\end{smallmatrix}\right) & i\text{ odd}\\[3mm]
\left(\begin{smallmatrix}\phantom{-}1&0\\-1&1\end{smallmatrix}\right)
& i\text{ even}
\end{array}\right.
$$
and let $\br_n$ act on homomorphisms $F_n\to\slz$ by Hurwitz automorphisms of
$F_n$.
Then the stabilizer group $Stab_{\psi_n}$ of
$\psi_n$ is equal to the stabilizer group $Stab_{\phi_n}$ of $\phi_n$.
\end{prop}

\proof
Both groups, $\slz$ and $\br_3$, are central extensions of $\psl$, and both
$\phi_n$ and $\psi_n$ induce the same homomorphism $\chi_n:F_n\to\psl$.
Of course $Stab_\chi$ contains $Stab_\phi$ and $Stab_\psi$ and thus our claim
is proved as soon as we can show the opposite inclusions.\\
First note that the braid action defined on homomorphisms as above is
equivalent to the Hurwitz action on the tuples of images of the specified
generators $t_i\in F_n$, hence the braid action will not change the conjugation
class of these images.\\
Now let $\b$ be a braid in $Stab_\chi$. Then
$\phi\circ\b(t_i)=(ab)^{3k_i}\phi(t_i)$ since $(ab)^3$ is the fundamental
element of $\br_3$ which generates the center of $\br_3$ and thus the
kernel of the extension $\br_3\to\psl$.
The degree homomorphism $d:\br_3\to\ZZ$ is a class function with value one on
all $\phi(t_i)$, hence $d((ab)^{3k_i})=0$. Since $d(ab)=2$ we conclude
$k_i=0$ and $\b\in Stab_\phi$.\\
Similarly we have $\psi\circ\b(t_i)=\pm\psi(t_j)$ for $\b\in Stab_\chi$.
Since the trace is a class function on $\slz$ which has value $2$ on all
$\psi(t_i)$ while it is $-2$ on $-\psi(t_i)$, we also get $\b\in Stab_\psi$.
\qed

\begin{prop}
\labell{slstab}
Let $F_n:=\langle t_i,1\leq i\leq n\,|\quad\rangle$ be the free group on
$n=l+l'$ generators, define a homomorphism $\psi_{l,l'}:F_n\to\slz$ by
$$
\psi_{l,l'}(t_i)=\left\{
\begin{array}{ll}
\left(\begin{smallmatrix}1&1\\0&1\end{smallmatrix}\right) & i>l, i\not\equiv l
\mod 2\\[2mm]
\left(\begin{smallmatrix}\phantom{-}1&0\\-1&1\end{smallmatrix}\right)
& i>l,i\equiv l \mod 2\\[2mm]
\left(\begin{smallmatrix}-1&\phantom{-}0\\\phantom{-}0&-1\end{smallmatrix}
\right) & i\leq l
\end{array}\right.
$$
and let $\br_n$ act on homomorphisms $F_n\to\slz$ by Hurwitz automorphisms of
$F_n$.
Then the stabilizer group $Stab_{\psi_{l,l'}}$ of $\psi_{l,l'}$
is generated by the image of $Stab_{\psi_{l'}}$ under the inclusion
$\br_{l'}\hookrightarrow\br_n$ mapping to braids with only the last $l'$
strands
braided and
$$
E_{l,l'}:=\left\langle \s_{ij}^{m_{ij}},\,1\leq i<j\leq n\,\left|\quad
m_{ij}=\left\{
\begin{array}{lcl}
1 & \text{if} & j\leq l\vee i\equiv j (2),i>l\\
2 & \text{ if } & i\leq l<j\\
3 & \text{ if } & i>l,i\not\equiv j(2)
\end{array}
\right.\right.\right\rangle.
$$
If $l'\leq6$ then even $Stab_{\psi_{l,l'}}=E_{l,l'}$.
\end{prop}

\proof
Again we argue with the equivalent Hurwitz action on images of the generators.
First we consider the induced action on conjugacy classes.
On $n$-tuples of conjugacy classes the Hurwitz action induces an action of
$\br_n$ through the natural homomorphism $\pi$ to the permutation group $S_n$.
Since the tuple induced from $\psi$ consists of $l$ copies of the conjugacy
class of $-id$ followed by $l'$ copies of the distinct conjugacy class of
$\psi(t_1)$, the associated stabilizer group is $\tilde E:=\pi\inv(S_{l}\times
S_{l'})$, and as in \cite{klui} one can check that
\begin{eqnarray*}
\tilde E & = & \langle \s_{ij}, i<j\leq l\text{ or }l<i<j;
\tau_{ij}:=\s_{ij}^2, i\leq l<j\rangle.
\end{eqnarray*}
So as a first step we have $Stab_\psi$ contained in $\tilde E$.\\

Since $-id$ is central it is the only element in its conjugacy class and we may
conclude that the $\tilde E$ orbit of $\psi$ contains only homomorphisms which
map the first $l$ generators onto $-id$. With a short calculation using that
$-id$ is a central involution we can check that the $\tau_{ij}$ act
trivially on
such elements:
$$
\begin{array}{cl}
& \tau_{ij}  (-id,...,-id,M_{l+1},...,M_{n})\\
= & \s_{i+1}\inv\cdots\s_{j-1}\inv\s_j^2\s_{j-1}\cdots\s_{i+1}
(-id,...,-id,M_{l+1},...,M_{n})\\
= & \s_{i+1}\inv\cdots\s_{j-1}\inv\s_j^2
(-id,...,-id,M_{l+1},...,M_{j-1},-id,M_j,...,M_{n})\\
= & \s_{i+1}\inv\cdots\s_{j-1}\inv
(-id,...,-id,M_{l+1},...,M_{j-1},-id,M_j,...,M_{n})\\
= &  (-id,...,-id,M_{l+1},...,M_{n})
\end{array}
$$
Therefore given $\b\in\tilde E$ as a word $w$ in the generators
$\s_{ij},\tau_{ij}$ of $\tilde E$ the action of $\b$ on $\psi$ is the same as
that of $\b'$ where $\b'$ is given by a word $w'$ obtained from $w$ by dropping
all letters $\tau_{ij}$.
By the commutation relations of the $\s_{ij}$ we may collect all letters
$\s_{ij},i,j\leq l$ to the right of letters $\s_{ij},i,j>l$ without changing
$\b'$ and get a factorization $\b'=\b_1'\b_2'$ with
$\b_1'\in\br_{l},\b_2'\in\br_{l'}$.\\

Hence $\b\in\tilde E$ acts trivially on $\psi$ if and only if $\b_1'\b_2'$ does
so if and only if $\b_2'$ acts trivially on $\psi_{l'}$.
Thus $Stab_{\psi_{l,l'}}$ is generated by the $\tau_{ij}$ the
$\s_{ij},i,j\leq l$
and the $\b_2'\in Stab_{\psi_{l'}}$.
Both conclusions of the proposition then follow since $\s_{ij},i,j>l$ are
contained in $Stab_{\psi_{l'}}$ and since they are even generators if
$l'\leq6$,
prop.\ \ref{l}.
\qed

\section*{mapping class groups of elliptic fibrations}

We return to elliptic fibrations and obtain some results concerning mapping
classes of elliptic fibrations. In fact we need to enrich the structure a
bit:

\begin{description}

\item[Definition:]
A {\it marked elliptic fibration} is an elliptic fibration with a distinguished
regular fibre, $f:X,F\to B,b_0$, which can be thought of as given by a {\it
marking} $F\hookrightarrow E$.

\item[Definition:]
A {\it fibration preserving map} of a marked elliptic surface $f:X,F\to B,b_0$
is a homeomorphism $\Phisl_X$ of $X$ such that $f\circ\Phisl_X=\bdiff_{B,b_0}\circ f$
for a homeomorphism $\bdiff_{B,b_0}$ of $(B,b_0)$ and such that $\Phisl_X|_F$ is
isotopic to the identity on $F$.\\
The map $\bdiff_{B,b_0}$ is called the {\it induced base homeomorphism}.

\end{description}

An induced homeomorphism necessarily preserves the set $\Delta(f)$ of singular
values of the fibration map $f$ and therefore can be regarded as a
homeomorphism of
the punctured base $B,\Delta(f)$ preserving the base point.\\

On the other hand with each elliptic fibration $f:X\to B$ we have a torus
bundle
over the complement $B^0$ of $\Delta(f)$. Its structure homomorphism is the
natural map
$$
\psi:\pi_1(B^0,b_0)\tto \diff(F)
$$
to the group of isotopy classes of diffeomorphisms of the distinguished fibre.

\begin{lemma}
\labell{moi-equiv}
Let $X,F\to B,b_0$ be a marked elliptic fibration and $\b$ a braid represented
by a homeomorphism of $B^0,b_0$. Then there is a fibration preserving map
$\Phisl_X$ inducing $\pi(\b)$ if and only if $\b$ stabilises the structure map
of the associated torus bundle.
\end{lemma}

\proof
A fibration preserving homeomorphism $\Phisl$ of an unmarked elliptic surface
induces
a map $\bdiff_B$ of the punctured base $B^0$. By the classification
of torus bundles there exists then a commutative diagram
\begin{eqnarray*}
\pi_1(B^0,b_0) & \stackrel{(\bdiff_B)_*}{-\!\!\!-\!\!\!\tto}
& \pi_1(B^0,\bdiff_B(b_0))\\
\downarrow\psi_{b_0} & & \downarrow\psi_{\bdiff(b_0)}\\
\diff(F) & \stackrel{(\Phisl|_F)_*}{-\!\!\!-\!\!\!-\!\!\!\tto}
& \diff(\Phisl(F))
\end{eqnarray*}
But the result of Moishezon \cite[p.\ 169]{moi} implies that the reverse
implication is true in the absence of multiple fibres.\\
If now $\Phisl$ is a fibration preserving homeomorphism of a marked elliptic
surface
then the bottom map is the identity and the claim is immediate.
\qed

In order to relate to the results of the last section we use surjective maps
$$
F_{6k+l}\to\pi_1(B,b_0)
$$
provided by a choice of geometric basis, i.e.\ an
ordered system of generators which are simultaneously represented by disjoint
loops, each going around a single element of $\Delta(f)$.

\begin{lemma}
\labell{moi-normal}
Given a marked elliptic fibration $X,F\to B,b_0$ with singular fibres only of
types $I_1,I^*_0$ and an isomorphism $\diff(F)\cong\slz$, there is a choice of
geometric basis for $\pi_1(B^0,b_0)$ such that the structure homomorphism of the
associated bundle is $\psi_{6k,l}$.
\end{lemma}

\proof
The proof proceeds along the lines of Moishezons proof \cite{moi}, cf.\
\cite{fm}, for the normal form of an elliptic surface with only fibres of type
$I_1$. The same strategy leads to our claim since fibres of type $I_0^*$ have
local monodromy in the center of $\slz$.
\qed

By now we have finally got all necessary results to prove theorem 1 as stated
in the introduction.\\

\proof[ of theorem \ref{1}]
As before $M^0_{6+l}$ denotes the mapping class group of
$B,\Delta(f)$.
We have previously shown that the mapping classes induced by fibration
preserving maps are represented by braids acting trivially on the structure
homomorphism of the torus bundle given with the elliptic fibration, lemma
\ref{moi-equiv}.

By lemma \ref{moi-normal} and prop.\ \ref{slstab} the corresponding group is
conjugation equivalent to $\pi(E_{6,l})=\overline{E}_{6,l}$.
On the other hand the monodromy group is in the conjugation class of
$\overline{E}_{6,l}$ by the main theorem. Moreover for each mapping
class of the monodromy group there is by prop.\
\ref{lift} a fibration preserving diffeomorphism, so we get an inclusion and
hence both groups coincide as claimed.
\qed

\end{document}